\documentclass[12pt]{iopart}
\eqnobysec
 
\begin{document}
\title{Inverse scattering transform for the Toda lattice with steplike initial data}
\author{Ag.Kh. Khanmamedov$^1$$^,$$^2$}
\address{$^1$ Institute of Mathematics and Mechanics of NAS of Azerbaijan, 9 F Agayev str., AZ1141, Baku, Azerbaijan }
\address{$^2$ Institute Applied Mathematics, Baku State University, 23 Z Khalilov str., AZ1148, Baku, Azerbaijan }
\ead{agil\_khanmamedov@yahoo.com}
\begin{abstract}
We study the solution of the Toda lattice Cauchy  problem with
steplike initial data. The initial data are supposed to tend to zero
as $n\rightarrow+\infty$. By the inverse scattering transform method
formulas allowing us to find solution of the Toda lattice is
obtained.

\end{abstract}
\ams{34K29; 35Q58}  \maketitle

\section{ Introduction}

The Toda lattice has some very important applications in the theory
of physics  of nonlinear processes (see [1]). It is known the
inverse scattering method allows one to investigate in detail the
Cauchy problem for the Toda lattice in the different classes of
initial data (see [1]-[15] and references therein).
 The last  problem for the doubly-infinite Toda lattice

\begin{equation} \label{GrindEQ__1_1_}
\left\{\begin{array}{l} {\dot{a}_{n} =\frac{a_{n} }{2} \left(b_{n+1}
-b_{n} \right),\, \, \, \, \, \cdot =\frac{d}{dt} ,a_{n} =a_{n}
(t)>0,} \\ {\dot{b}_{n} =a_{n}^{2} -a_{n-1}^{2} ,\, \, \, \, \, \,
b_{n} =b_{n} (t),\, \, \, n=0,\pm 1,\pm 2,....\, \, }
\end{array}\right.
\end{equation}
with fast stabilized or steplike fast stabilized initial data is
investigated in [1]-[9] (see also references therein) by the method
of inverse scattering transform. However, this problem is not
studied in the case of steplike initial data, where $a_n$ tend to
zero as $n\to +\infty $ (or   $n\to -\infty $).

 In this paper we study the Cauchy problem for the system (1.1) with initial data

\[a_{n} (0)\to 0,\, \, b_{n} (0)\to 0\quad as\, \, \, \, \, n\to +\infty ,\]

\begin{equation} \label{GrindEQ__1_2_}
\sum _{n<0}^{} \left|n\right|\left\{\left|a_{n}
(0)-1\right|+\left|b_{n} (0)\right|\right\}<\infty \, \, .
\end{equation}
The solution is considered in the class
\[\|a_n(t)\|_{C[0,T]}\rightarrow 0,\,\,\|b_n(t)\|_{C[0,T]}\rightarrow
0,\,\,\,\texttt{ as}\,\,\, n\rightarrow +\infty,\]
\begin{equation} \label{GrindEQ__1_3_}
\left\| Q(t)\right\| _{C[0,\, \, T]} <\infty \, \,,
\end{equation}
for arbitrary $T>0$, where

\[Q(t)=\sum _{n<0}^{} \left|n\right|\left(\left|a_{n} (t)-1\right|+\left|b_{n} (t)\right|\right)\, \, .\]

 Note, we cannot apply directly method given in [1]-[9] for the case $\mathop{\inf }\limits_{n} a_{n} >0$, because
  the Jost solution with the asymptotic behaviour on an $+\infty $ does not exist in our case. On the other hand,
   method of inverse problem   is used (see [10]) in the case when Jacobi operator associated with (1.1) has the continuous spectrum $\left[a,\, \, b\right]$ of multiplicity two. But this method cannot be used when the spectrum of the Jacobi operator has a continuous spectrum of multiplicity one and a discrete spectrum.

 The paper is organized as follows. In section 2 we formulate some auxiliary facts to the inverse scattering problem
   for the Jacobi operator associated with (1.1)-(1.2). In section 3 we describe the evolution of the scattering data of problem (1.1)-(1.2) .

In the last section we prove existence of the solution of the
problem (1.1)-(1.2) in class (1.3).

\section{ The scattering problem}

 Consider Jacobi operator $L$ generated in $\, \ell ^{2} \left(-\infty ,\, \infty \right)$ by the finite-difference operations

\[\left(Ly\right)_{n} =a_{n-1} y_{n-1} +b_{n} y_{n} +b_{n} y_{n+1} ,\]
in which the real coefficients $a_{n} >0,\; b_{n} $ satisfy the
conditions

\[a_{n} \to 0,\, \, b_{n} \to 0\quad as\, \, \, \, \, \, n\to +\infty ,\]

\[\sum _{n<0}^{} \left|n\right|\left\{\left|a_{n} -1\right|+\left|b_{n} \right|\right\}<\infty \, \, .\]
The interval $[-2, \, 2]$ is the  continuous spectrum of
multiplicity one of operator $L$ (see [16],[17]). Beyond the
continuous spectrum, $L$ can have a finite number of simple
eigenvalues $ \mu_k(t),\,\, k=1,...,p$.

 Let us formulate some auxiliary facts related to the inverse scattering problem for the equation

\begin{equation} \label{GrindEQ__2_1_}
\left(Ly\right)_{n} =\lambda y_{n} ,\quad n=0,\, \, \pm 1,...,\, \,
\, \lambda \in \mathbf{C}
\end{equation}
Many of these facts can be found in [16],[17].

 Let $P_{n} (\lambda )$ and $Q_{n} (\lambda )$ be solutions
of  Eq. (2.1)  with initial conditions

\[P_{-1} (\lambda )=0, \,\, P_{0} (\lambda )=1,\]

\[Q_{0} (\lambda )=0, \, Q_{1} (\lambda )=\frac{1}{a_{0}^{} } .\]
We denote by $L_{0} $ semi-infinite Jacobi operator generated $\,
\ell ^{2} \left[0,\, \, \infty \right)$ by Eq. (2.1) as $n\ge 0$ and
the boundary condition   $y_{-1}^{} =0$. This operator is completely
continuous.  Moreover, the spectral function $\rho (\lambda )$ of
$L_{0} $ represented [18] in the form

\[\rho (\lambda )=\sum _{\lambda _{n} <\lambda }\beta _{n}^{-2}  \]
where $\lambda _{n} $ is the eigenvalue of $L_{0} $ and $\beta _{n}
 $ is the norm of the eigenfunction corresponding to the
$\lambda _{n}$.

 As is known from [18]-[19], the right Weyl function of the problem (2.1) has the form

\begin{equation} \label{GrindEQ__2_2_}
m(\lambda )=\int _{-\infty }^{\infty }\frac{d\rho (\tau )}{\tau
-\lambda }  ,
\end{equation}
or

\[m(\lambda )=\sum _{n=1}^{\infty } \frac{\beta {}_{n}^{-2} }{\lambda _{n} -\lambda } .\]
 where  $\lambda _{n} \to 0$ as  $n\to \infty $. It follows from [12]-[13] that for
  $\lambda \ne \lambda _{k} ,\, \, \, \, k=1,2,...,$ Eq.(2.1) has Weyl solution

\begin{equation} \label{GrindEQ__2_3_}
\psi _{n} (\lambda )=\, \, Q_{n} (\lambda )+m(\lambda )P_{n}
(\lambda ),
\end{equation}
``on the right semiaxis'' (such that $\sum\limits _{n=0}^{\infty }
\left|\psi _{n} \left(\lambda \right)\right|^{2} <\infty .$ ).

 Suppose that   $\Gamma $\, is the complex $\lambda $--plane with cut along the interval $[-2,\,2]$. In the  plane   $\Gamma $,   consider the function

\[z(\lambda )=\frac{\lambda }{2} +\sqrt{\frac{\lambda ^{2} }{4} -1} \]
choosing the regular branch of the radical so that
$\sqrt{\frac{\lambda ^{2} }{4} -1} <0$ for $\lambda >2$. We often
omit the dependence of $z(\lambda )\, \, $on $\, \lambda $ in what
follows. Thus, in the formulas involving $z$ and $\, \lambda $, we
always assume that $z$ is as in the above equation.

 It is well known (see, for example, [20]) that Eq.(2.1) has a Jost solution represented in
 the form

\begin{equation} \label{GrindEQ__2_4_}
f_{n} (\lambda )=\alpha _{n} z^{-n} \left(1+\sum _{m<0}^{} A_{nm}
z^{-m} \right).
\end{equation}
 The coefficients are given by

\begin{equation} \label{GrindEQ__2_5_}
a_{n} =\frac{\alpha _{n} }{\alpha _{n+1} } ,\, \, \, \, b_{n}
=A_{n,-1} -A_{n+1,-1}.
\end{equation}

 Without restriction of generality we can suppose that $\lambda _{m} \in \left(-2,\, \, \, 2\right)$ for
 any $m=1,2,....$. As known [16],[17], for  $\lambda \in \partial \Gamma ,\, \, \, \lambda ^{2} \ne 4,\, \, \,
 \lambda \ne \lambda _{m} $ identity

\begin{equation} \label{GrindEQ__2_6_}
\psi _{n} (\lambda )=a(\lambda )\overline{f_{n} (\lambda
)}+\overline{a(\lambda )}f_{n} (\lambda )
\end{equation}
holds,    where  the function $a(\lambda )$ can be regularly
continued  to $\Gamma $. Note also, $a(\lambda )$ can have a finite
number of coinciding simple zeros outside the interval $\left[-2,\,
\, \, 2\right]$, because, these zeros constitute the discrete
spectrum $\mu _{k} ,\, \, \, k=1,...,p,$ of the operator $L$.

 Introduce  reflection $R(\lambda )$ coefficient by the formula

\[ R(\lambda )=\frac{\overline{a(\lambda )}}{a(\lambda )} .\]
The function $R(\lambda$ is continuous for
$\lambda\in\partial\Gamma.$ Setting $n=-1$ and $n=0$ in the identity
(2.6) yields the expression
\begin{equation} \label{GrindEQ__2_7_}
m(\lambda)=-\frac{1}{a_{-1}}
\frac{\overline{f_0(\lambda)}+R(\lambda)f_0(\lambda)}{f_{-1}(\lambda)+R(\lambda)f_{-1}(\lambda)}
\end{equation}
 The norming constants $M_{k} (t)$
corresponding to the $\mu_k(t)$ are given as

\[M_{k}^{-2} =\sum _{n=-\infty }^{\infty }f_{n}^{2}  (\mu _{k} )\, \, ,\, \, k=1,...,p.\]

 The set of quantities $\left\{R(\lambda );\mu _{k} ;\, \, M_{k} ,\, \, k=1,....,p\, \right\}$ is called the scattering data for the Jacobi operator $L$. The inverse scattering problem for $L$ is to recover the coefficients $a_{n} ,\; b_{n} $ from the scattering data.

In solving the inverse problem, an important role is played by the
Marchenko-type basic equation. Define

\begin{equation} \label{GrindEQ__2_8_}
F_{n} =\sum _{k=1}^{p}M_{k}^{-2} z_{k}^{-n}  +\frac{1}{2\pi i} \int
_{\partial \Gamma } \frac{R(\lambda )}{z^{-1} -z}  z_{}^{-n}
d\lambda ,
\end{equation}
where $z_{k} =z(\mu _{k} ), \, \, k=1,...,p$.

 Then $A_{nm} $ and $\alpha _{n} $ involved in (2.4) satisfy the relations

\begin{equation} \label{GrindEQ__2_9_}
F_{2n+m} +A_{nm} +\sum _{k<0}^{} A_{nk} F_{2n+m+k} =0,\, \, \, \, \,
m<n\le 0\, \, ,
\end{equation}

\begin{equation} \label{GrindEQ__2_10_}
\alpha _{n}^{-2} =1+F_{2n} +\sum _{k<0}^{} A_{nk} F_{2n+k} ,\, \, \,
\, \, \, n\le 0\, \, .
\end{equation}
To reconstruct the operator $L$, we consider Eq.(2.8) which is
constructed by the scattering data. We find $A_{nm} $ and $\alpha
_{n} $ from Eqs.(2.9) and (2.10), respectively, the firs one having
a unique solution with respect to $A_{nm} $.The coefficients $a_{n}
\, \, \, \, and\, \, \, \; b_{n} $ are definides for $n<0$ by (2.5).
$f_{n} (\lambda )$ for $n\le 0$ are defined by (2.4). From the
formula (2.7) we obtain Weyl function $m(\lambda )$. The spectral
measure $d\rho (\lambda )$ can be found by the formula

\[d\rho (\lambda _{n} )=\mathop{\lim }\limits_{\lambda \to \lambda _{n} }(\lambda_n-\lambda )m(\lambda),  \,\,\, n=1,2,...
.\] Using the approach in [12],[13],[19], we can reconstruct
semi-infinite Jacobi operator $L_{0} $  by its spectral measure
$d\rho (\lambda )$. Therefore, we find  $a_{n} ,\; b_{n} $ for $n\ge
0$.

\section{ Evolution of the scattering data}

In this section we use the inverse scattering transform method to
solve the problem (1.1)-(1.2). Let $a_{n} (t),\; b_{n} (t)$ be a
solution of the problem (1.1)-(1.2) satisfying (1.3). Consider the
Jacobi operator $L=L(t)$ associated with  $a_{n} =a_{n} (t),\; b_{n}
=b_{n} (t)$. Jost and Weyl solutions, reflection coefficient,
spectral measure now depend on the additional parameter $t\in [0,\,
\, \, \infty )$.

 \textbf{  Theorem 1.  }\emph{If the coefficients $a_{n} =a_{n} (t),\; b_{n} =b_{n}
  (t)$ of Eq.(2.1) are solutions to problem (1.1)-(1.2) in the class (1.3), then the evolution of the scattering data is  described by the
 formulas}

\begin{equation} \label{GrindEQ__3_1_}
R(\lambda ,t)=R(\lambda ,0)\, e^{\left(z^{-1} -z\right)t}
\end{equation}

\begin{equation} \label{GrindEQ__3_2_}
\mu_k(t)=\mu_k(0),\,\, k=1,...,p
\end{equation}

\begin{equation} \label{GrindEQ__3_3_}
M_{k}^{-2} (t)=M_{k}^{-2} (0)e^{\left(z^{-1}_k -z_k\right)t} ,\, \,
\, \, \, \, \, \, \, z_{k} =z(\mu _{k} )\, \, ,\, \, k=1,...,p.
\end{equation}

\textbf{Proof.}  System (1.1) is represented (see, for example
[8],[13]) in the Lax form

\begin{equation} \label{GrindEQ__3_4_}
\dot{L}=\left[L,\, \, A\right]=AL-LA,
\end{equation}
where  $A=A(t)$ are Jacobi operator in $\, \ell ^{2} \left(-\infty
,\, \infty \right)$:

\[\left(Ay\right)_{n} =\frac{1}{2} a_{n} y_{n+1} -\frac{1}{2} a_{n-1} y_{n-1} .\]
Since (3.4) implies that the family of operators $L=L(t)$ are
unitarily equivalent (see [5],[8]), the spectrum of $L=L(t)$ does
not depend on  $n$ and (3.2) is valid.

 Let $f_{n} (\lambda ,\, \, t)$and $\psi _{n} (\lambda ,\,
\, t)$respectively be the Jost and  Weyl solutions of the Eq.(2.1)
with the parameter $t$. Consider the indetity (2.6) with the
parameter $t$. As follows from [8],[12] the function $\frac{d}{dt}
\psi _{n} -\left(A\psi \right)_{n} $ is also a solution of the
Eq.(2.1) with the parameter $t$. Appling the operator $\frac{d}{dt}
-A$ to (2.6), taking into account that the Jost solution $f_{n}
(\lambda ,\, \, t)$ does not depend (see [8], on $t$ asymptotically,
we obtain

\begin{equation} \label{GrindEQ__3_5_}
\begin{array}{l} {\frac{d}{dt} \psi _{n} -\left(A\psi \right)_{n} =\left(\dot{a}(\lambda ,\, \, t)
+\frac{1}{2} \left(z^{-1} -z\right)a(\lambda ,\, \, t)\right)\overline{f_{n} (\lambda ,\, \, t)}+} \\
\\ {+\left(\overline{\dot{a}(\lambda ,\, \, t)}-\frac{1}{2}
\left(z^{-1} -z\right)\overline{a(\lambda ,\, \, t)}\right)f_{n}
(\lambda ,\, \, t).} \end{array}
\end{equation}
On the other hand, we find

\[\frac{d}{dt} P_{0} -\left(AP\right)_{0} =\frac{b_{0} -\lambda }{2} ,\, \, \, \frac{d}{dt} P_{-1} -\left(AP\right)_{-1} =-a_{-1} \, \, .\]
Since $P_{n} (\lambda ,\, \, t)$ and $Q_{n} (\lambda ,\, \, t)$ are
linearly independent, the function $\frac{d}{dt} P_{n}
-\left(AP\right)_{n} $ can be represented as

\[\frac{d}{dt} P_{n} -\left(AP\right)_{n} =A(\lambda ,\, \, t)P_{n} +D(\lambda ,\, \, t)Q_{n} .\]
Setting $n=-1$ and $n=0$ in the last relation, we find that

\[A(\lambda ,\, \, t)=\frac{b_{0} -\lambda }{2} ,\, \, \, \, D(\lambda ,\, \, t)=a_{-1}^{2} .\]
Therefore,

\[\frac{d}{dt} P_{n} -\left(AP\right)_{n} =\frac{b_{0} -\lambda }{2} P_{n} +a_{-1}^{2} Q_{n} .\]
The same arguments are valid for solution $Q_{n} (\lambda ,\, \,
t).$
 Thus, we have the formula

\[\frac{d}{dt} Q_{n} -\left(AQ\right)_{n} =-P_{n} +\frac{\lambda -b_{0} }{2a_{0} } Q_{n} .\]
Now by the formula (2.3) with the parameter $t$ we find that

\begin{equation} \label{GrindEQ__3_6_}
\begin{array}{l} {\frac{d}{dt} \psi _{n}-\left(A\psi \right)_{n} =\left(a_{-1}^{2} m(\lambda ,\, \, t)
+\frac{\lambda -b_{0} }{2a_{0} } \right)Q_{n} +} \\
\\ {+\left(\dot{m}(\lambda ,\, \, t)+\frac{b_{0}
 -\lambda }{2}m(\lambda,\,t) -1\right)P_{n} .} \end{array}
\end{equation}
Since $L=L(t)$  is selfadjoint and bounded,  $\frac{d}{dt} \psi _{n}
-\left(A\psi \right)_{n} $ must satisfy the relation

\begin{equation} \label{GrindEQ__3_7_}
\frac{d}{dt} \psi _{n} -\left(A\psi \right)_{n} =\theta (\lambda ,\,
\, t)\psi _{n}.
\end{equation}
Hence, we can represent the function $\frac{d}{dt} \psi _{n}
-\left(A\psi \right)_{n} $ as

\begin{equation} \label{GrindEQ__3_8_}
\frac{d}{dt} \psi _{n} -\left(A\psi \right)_{n} =\theta (\lambda ,\,
\, t)Q_{n} +\theta (\lambda ,\, \, t)m(\lambda ,\, \, t)P_{n} .
\end{equation}
Comparing this identity with (3.6), we have

\begin{equation} \label{GrindEQ__3_9_}
\theta (\lambda ,\, \, t)=a_{-1}^{2} m(\lambda ,\, \,
t)+\frac{\lambda -b_{0} }{2a_{0} } ,
\end{equation}

Further, according to (2.6), (3.5), (3.7),

\[\begin{array}{l} {\theta (\lambda ,\, \, t)a(\lambda ,\, \, t)\overline{f_{n} }+\theta (\lambda ,\, \, t)
\overline{a(\lambda ,\, \, t)}f_{n} =\left(\dot{a}(\lambda ,\, \,
t)+\frac{1}{2} \left(z^{-1} -z\right)a
(\lambda ,\, \, t)\right)\overline{f_{n} }+} \\
\\{+\left(\overline{\dot{a}(\lambda ,\, \, t)}-\frac{1}{2}
\left(z^{-1} -z\right)\overline{a(\lambda ,\, \, t)}\right)f_{n} .}
\end{array}\] Since $f_{n} $ and $\overline{f_{n} }$ are linearly
independent, so substituting (3.9) into the last indetity, we obtain

\[\begin{array}{l} {\dot{a}(\lambda ,\, \, t)+\frac{1}{2} \left(z^{-1} -z\right)a(\lambda ,\, \, t)=
\left(a_{-1}^{2} m(\lambda ,\, \, t)+\frac{\lambda -b_{0} }{2a_{0} } \right)a(\lambda ,\, \, t),} \\
\\ {\overline{\dot{a}(\lambda ,\, \, t)}-\frac{1}{2} \left(z^{-1}
-z\right)\overline{a(\lambda ,\, \, t)}=\left(a_{-1}^{2} m(\lambda
,\, \, t)+\frac{\lambda -b_{0} }{2a_{0} } \right)\overline{a(\lambda
,\, \, t)}.} \end{array}\] From this relations, we get

\[\dot{R}(\lambda,\,t)=(z^{-1}-z)R(\lambda,\,t),\]
 which imply
(3.1).

 Now, let $g_{n} (\mu _{k} ,\, \, t)$ be a normalized eigenfunction of $L$. Since the eigenvalues $\mu_k, k=1,...,p$,
 of this operator are simple, we have

\[\frac{d}{dt} g_{n} -\left(Ag\right)_{n} =cg_{n} .\]
Taking the scalar products of  $g_{n} $ with both sides of this
equality in
  $\, \ell ^{2} \left(-\infty ,\, \infty \right)$ and using
    $\, \left\| \psi _{n} \right\| _{\ell ^{2} \left(-\infty ,\, \infty \right)} =1$ and $A^{*} =-A,$
    we obtain $c=0$. Therefore,

\begin{equation} \label{GrindEQ__3_10_}
\frac{d}{dt} g_{n} -\left(Ag\right)_{n} =0
\end{equation}
On the other hand, if a normalized eigenfunction $g_{n} (\mu _{k}
,\, \, t)$ corresponds to the eigenvalue $\mu _{k} $, then

\[g_{n} (\mu _{k} ,\, \, t)=c_{k} (\, t)f_{n} (\mu _{k} ,\, \, t).\]
This implies that $M_{k}^{2} (t)=c_{k}^{2} (t)$. By virtue of (2.4),
we find that

\[\frac{d}{dt} g_{n} -\left(Ag\right)_{n} \sim \left(\dot{c}_{k} (t)+\frac{z_{k} -z_{k}^{-1} }{2} c_{k} (t)\right)z_{k}^{-n} \]
as   $n\to -\infty $.
  Taking into account (3.10), we have

\[\dot{c}_{k} (t)+\frac{z_{k} -z_{k}^{-1} }{2} c_{k} (t)=0\]
This equation implies the relation (3.3).

 The theorem is proved.

Using Theorem 1, we obtain the following procedure for solving
problem (1.1),(1.2) based on the inverse scattering transform
method: Initial data (1.2) is given. Construct $R(\lambda
,0)$,$\mu_k(0)$, $M_{k} (0),\, \, k=1,....,p$ . Calcullate
$R(\lambda ,t)$,$\mu_k(t)$,$M_{k} (t)$ using formulas (3.1)-(3.3).
Construct a solution by solving the inverse problem by applying
approach of the section  2 with $R(\lambda ,0)$,$\mu_k(0)$, $M_{k}
(0),\, \, k=1,....,p$ replaced by (3.1)-(3.3).

\section{ Solvability of the Cauchy problem for the Toda lattice}

In section 3, while constructing a solution to problem (1.1)-(1.2),
we assumed that this solution exists in the class (1.3). Let us now
investigate its existence.

\textbf{Theorem 2}. \emph{The problem (1.1)-(1.2) has a unique
solution in the class (1.3).}

\textbf{Proof.} Denote by $B_{} $ the Banach space of pairs of
sequences $y=\left(y_{1,n} ,\, \, y_{2,n} \right)_{ n=-\infty
}^{\infty } $ for which the norm $\left\| y\right\| _{B}
=\mathop{\sup\limits_{n\ge 0} \left(\left|y_{1,n}
\right|+\left|y_{2,n} \right|\right)+\sum\limits _{n<0}
\left|n\right|\left(\left|y_{1,n} \right|+\left|y_{2,n}
\right|\right)}\, $ is finite. Then (see [21]) the set
$C\left(\left[0,\, \, T\right];\, \, B\right)$ of the continuous on
an interval $\left[0,\, \, T\right]$ with respect to the norm
$\left\| \cdot \right\| _{B} $ functions is the Banach space.

 Let as assume that

\begin{equation} \label{GrindEQ__4_1_}
\begin{array}{l} {x_{1,n} =\left\{\begin{array}{l} {a_{n} (t)\, \, \, \, for\, \, \, \, \, \, \, n\ge 0,} \\
\\ {a_{n} (t)-1\, \, \, \, for\, \, \, n<0,\, } \end{array}\right. }
\\
\\ {x_{2,n} =b_{n} (t).} \end{array}
\end{equation}
Then system (1.1) is equivalent to the system

\begin{equation} \label{GrindEQ__4_2_}
\left\{ \begin{array}{l}
 \dot x_{1,n}  = \frac{1}{2}x_{1,n} \left( {x_{2,n + 1}  - x_{2,n} }
  \right) + \frac{1}{2}\left( {1 - \delta _{n,\left| n \right|} } \right)\left( {x_{2,n + 1}  - x_{2,n} } \right), \\
  \\
 \dot x_{2,n}  = x_{1,n}^2  - x_{1,n - 1}^2  + 2\left( {1 - \delta _{n,\left| n \right|} } \right)\left( {x_{1,n}
  - x_{2,n - 1} } \right),
 \end{array} \right.
\end{equation}

where $\delta _{n,m} $ is the Kronecker symbol.

 Denote by ${F}$ the operator generated the right-hand sides of system (4.2). Note, operator  ${F}$
 is strongly continuously differentable in the space  $C\left(\left[0,\, \, T\right];\, \, B\right)$.

 Now passing to the integral equation in the standard manner, we find problem (4.2) with initial conditions

\begin{equation} \label{GrindEQ__4_3_}
\begin{array}{l} {x_{1,n} (0)=\left\{\begin{array}{l} {a_{n} (0)\, \, \, \, for\, \, \, \, \, \, \, n\ge 0,} \\
\\ {a_{n} (0)-1\, \, \, \, for\, \, \, n<0,\, } \end{array}\right. }
\\
\\ {x_{2,n} (0)=b_{n} (0).} \end{array}
\end{equation}
 is equivalent to the equation

\begin{equation} \label{GrindEQ__4_4_}
x(t)=x(0)+\int _{0}^{t}{F}(x(\tau ))d\tau
\end{equation}
Applying the principe of compressed maps, we find that problem (4.4)
on some interval $\left[0,\, \, \delta \right]$ has a unique
solution $x(t)$with finite norm $\left\| x(t)\right\|
_{C\left(\left[0,\, \, \delta\right];\, \, B\right)} <\infty $. Let
us show that this solution can be extended to the entire positive
semiaxsis. Assume the opposite. Then there exists a point $t^{*} \in
\left(0,\, \, \infty \right)$such that problem (4.2)-(4.3) has a
solution $x(t)=\left(x_{1,n} (t),\, \, x_{2,n} (t)\right)$on the
interval $\left[0,\, \, t^{*} \right)$ but  $\overline{\mathop{\lim
}\limits_{t\to t^{*} -0} }\left\| x(t)\right\| _{B} =\infty $. It
follows from [8],[13] problem (1.1)-(1.2) has a unique
solution$\left(a_{n} (t),\, \, b_{n} (t)\right)$ in $C^{\infty }
\left(\left[0,\, \, \infty \right);\, \, M\right)$, where $\, M=\ell
^{\infty } \left(-\infty ,\, \infty \right)\oplus \ell ^{\infty }
\left(-\infty ,\, \infty \right)$. Hence, according to the (4.1)
problem (4.2)-(4.3) has a unique solution $x(t)=\left(x_{1,n} (t),\,
\, x_{2,n} (t)\right)$ satisfying

\[\left|x_{1,n} (t)\right|+\left|\, x_{2,n} (t)\right|<C\]
for  any $t\in \left[0,\, \, \infty \right)$, where $C$ does not
depend on $t$ . We integrate the system (4.2) over a interval
$\left[0,\, \, t\right]$. Then, using the last inequality, after
some simple transformations , we get

\[\left\| x(t)\right\| _{B} \le 2\left\| x(0)\right\| _{B} +\left(4C+4\right)\int _{0}^{t}\left\| x(\tau )\right\| _{B} d\tau  ,\, \, \, 0<t<t^{*} \, \, ,\]
which, according to the Gronwall's inequality implies

\[\left\| x(t)\right\| _{B} \le 2\left\| x(0)\right\| _{B}e^{\left(4C+4\right)t} .\]
Therefore, our assumption that $\overline{\mathop{\lim
}\limits_{t\to t^{*} -0} }\left\| x(t)\right\| _{B} =\infty $ is not
correct and problem (4.2)-(4.3) has a unique solution
$x(t)=\left(x_{1,n} (t),\, \, x_{2,n} (t)\right)\in
C\left(\left[0,\, \, T\right];\, \, B\right)$ for any $T>0$.
Integrating the system (1.1) over a interval $[0,\,t]$ and using
(4.1), we obtain that problem (1.1)-(1.2) be uniquely solvable in
the class (1.3).

 Thus, the theorem is proved.

\end{document}